\documentclass [a4paper,11pt]{amsart}

\usepackage[utf8]{inputenc}
\usepackage{xcolor}
\usepackage[pdftex]{graphicx} 
\usepackage{hyperref}
\DeclareGraphicsExtensions{.png,.pdf,.jpg}
\usepackage{amsmath}
\usepackage{amssymb}
\usepackage{amsxtra}
\usepackage{latexsym}
\usepackage[all]{xy}
\usepackage{array}
\usepackage{bbm}
\usepackage{tikz-cd}
\usepackage{stmaryrd}
\usepackage{float}
\RequirePackage{filecontents}    
\def\Ker{\operatorname{Ker}}

\def\dim{\operatorname{dim}}

\newcommand{\C}{\mathbb{C}}
\newcommand{\Q}{\mathbb{Q}}
\newcommand{\PP}{\mathbb{P}}
\newcommand{\R}{\mathbb{R}}
\newcommand{\K}{\mathbb{K}}
\newcommand{\Z}{\mathbb{Z}}

\newcommand{\hH}{\mathcal{H}}

\theoremstyle{definition}
\newtheorem{teo}{Theorem}[section]
\newtheorem{defi}[teo]{Definition}

\newtheorem{lem}[teo]{Lemma}
\newtheorem{propo}[teo]{Proposition}

\theoremstyle{definition}

\theoremstyle{definition}
\newtheorem{re}[teo]{Remark}

\begin{document}

\title[The Brasselet-Sch\"urmann-Yokura conjecture]{The Brasselet-Sch\"urmann-Yokura conjecture on $\boldsymbol L$-classes of projective varieties}

\author{J. Fern\'andez de Bobadilla}
\address{Javier Fern\'andez de Bobadilla:  
(1) IKERBASQUE, Basque Foundation for Science, Maria Diaz de Haro 3, 48013, 
    Bilbao, Basque Country, Spain;
(2) BCAM,  Basque Center for Applied Mathematics, Mazarredo 14, 48009 Bilbao, 
Basque Country, Spain; 
(3) Academic Colaborator at UPV/EHU.} 
\email{jbobadilla@bcamath.org}

\author{I. Pallar\'es}
\address{Irma Pallar\'es: KU Leuven, Celestijnenlaan 200B, B-3001 Leuven, Belgium.} 
\email{irma.pallarestorres@kuleuven.be}


\thanks{J.F.B. was supported by ERCEA 615655 NMST Consolidator Grant, MINECO by the project 
reference MTM2016-76868-C2-1-P (UCM), by the Basque Government through the BERC 2018-2021 program and Gobierno Vasco Grant IT1094-16, by the Spanish Ministry of Science, Innovation and Universities: BCAM Severo Ochoa accreditation SEV-2017-0718, and by VIASM, through a funded research visit.}
\thanks{I.P. was supported by ERCEA 615655 NMST Consolidator Grant, by the Basque Government through the BERC 2018-2021 program and Gobierno Vasco Grant IT1094-16, by the Spanish Ministry of Science, Innovation and Universities: BCAM Severo Ochoa accreditation SEV-2017-0718, by 12B1423N from FWO, Research Foundation Flanders.
}



\begin{abstract}
In 2010, Brasselet, Schürmann and Yokura conjectured an equality of characteristic classes of singular varieties between the Goresky-MacPherson $L$-class $L_*(X)$ and the Hirzebruch homology class $T_{1,*}(X)$ for a compact complex algebraic variety $X$  that is a rational homology manifold. In this note we give a proof of this conjecture for projective varieties based on cubical hyperresolutions, the Decomposition Theorem, and Hodge theory. The crucial step of the proof is a new characterization of rational homology manifolds in terms of cubical hyperresolutions which we find of independent interest. 
\end{abstract}

\maketitle


\section{Introduction}\label{intro}


\subsection{Characteristic classes, from manifolds to singular varieties} In \cite{H:1966}, F. Hirzebruch defined a theory of characteristic classes of vector bundles depending on a parameter $y$ unifying Chern classes, Todd classes and Thom-Hirzebruch $L$-classes. The {\em Hirzebruch (cohomology) class} $T_y^*$ of a complex manifold $X$ is formally given by: 
 $$T_y^*(X)=\prod_{i=1}^{\dim X} \frac{\alpha_i(1+y)}{1-e^{-\alpha_i(1+y)}-\alpha_i y}\in\Q[y][[\alpha_1,\dots,\alpha_{\dim X}]]$$ where $\alpha_i$ are the Chern roots of the tangent
bundle $TX$ of $X$. It specializes to the Chern class of $X$ for $y=-1$, to the Todd class of $X$ for $y=0$, and to the Thom-Hirzebruch $L$-class of $X$ for $y=1$. F. Hirzebruch also defined the  $\chi_y$-{\em characteristic} of a compact complex manifold $X$ as$$ \chi_y(X)=\sum_{p,q\in \mathbb N}(-1)^q \dim H^q(X,\Omega_X^p) y^p \in \mathbb Z[y],$$ and proved a {\em generalized Hirzebruch-Riemann-Roch theorem} which states
$$\chi_y(X)=\int_X T_y^*(X)\cap [X].$$
This theorem specializes to the Gauss-Bonnet-Chern theorem for $y=-1$, to the Hirzebruch-Riemann-Roch theorem for $y=0$, and to the Hirzebruch signature theorem for $y=1$.

\medskip

Chern classes, Todd classes and Thom-Hirzebruch $L$-classes of manifolds were individually extended to singular varieties.  The first extension of Chern classes for singular varietes was introduced by M. H. Schwartz in \cite{Schw:1965}, using obstruction theory. In \cite{McP:1974}, R. MacPherson, answering a conjecture of P. Deligne and A. Grothendieck, proved the unique existence of a natural transformation 
$$c_*\colon F(-)\to H_*(-),$$ 
where $F(-)$ is the covariant functor of constructible functions  and $H_*(-)$ the usual Borel-Moore covariant $\mathbb Z$-homology functor, satisfying $c_*(1_X)=c^*(TX)\cap [X]$ if $X$ is a smooth variety, where $1_X$ is the constant function with value one on $X$. For any possibly singular variety $X$,  J.-P. Brasselet and M. H. Schwartz in \cite{BS:1981} proved that MacPherson's Chern class $c_*(1_X)$ is the same as the previously defined Schwartz's Chern class of $X$ via the Alexander isomorphism. The {\em Chern-Schwartz-MacPherson} class of any variety $X$ is defined as $c_*^{SM}(X):=c_*(1_X)$. 

In \cite{BFM:1975}, P. Baum, W. Fulton and R. MacPherson defined a Todd class for singular varieties by proving the unique existence of a natural transformation 
$$td_*\colon G_0(-)\to H_*(-)\otimes \mathbb Q,$$ where $G_0(-)$ is the covariant Grothendieck functor of coherent sheaves, satisfying
$td_*(\mathcal O_X)=td^*(TX)\cap [X]$,  if $X$ is a smooth variety, where $\mathcal O_X$ is the structure sheaf of $X$. The {\em Baum-Fulton-MacPherson Todd class} of any variety $X$ is defined as $td_*^{BFM}(X):=td_*(\mathcal O_X)$.

In \cite{GM1:1980}, \cite{GM2:1983}, M. Goresky and R. MacPherson introduced Intersection Homology obtaining a Poincaré duality and a definition of signature for singular varieties. Furthermore, the same authors defined an $L$-class for compact singular varieties using cohomotopy. S. E. Cappell, J. L. Shaneson and S. Weinberger in  \cite{CSW:1991}, \cite{CS:1991}, \cite{CS2:1991} (see also \cite{BSY:2010}) extended the methods of M. Goresky and R. MacPherson on the construction of an $L$-class for singular varieties.  They defined an $L$-class for self-dual complexes of sheaves by proving the unique existence of a natural transformation 
$$L_*\colon \Omega_{\K}(-)\to H_*(-;\mathbb Q),$$ 
where $\Omega_{\K}(-)$ is the covariant cobordism functor of self-dual $\K$-complexes of sheaves where $\K$ is a subfield of $\R$, satisfying  $L_*([\K_X[\dim_{\mathbb C}X]])=L^*(TX)\cap [X]$  if $X$ is  a compact smooth variety. The {\em Cappell-Shaneson-Weinberger $L$-class} of a compact complex algebraic variety $X$ is defined as $L_*(X):=L_*([IC_X])$, where $IC_X$ denotes the intersection cohomology sheaf complex on $X$, which coincides with the previously defined Goresky-MacPherson $L$-class, see also \cite{Yo:1995}.

\subsection{Brasselet-Schürmann-Yokura's Hirzebruch classes of singular varieties} In \cite{BSY:2010}, J.-P. Brasselet, J. Schürmann and S. Yokura unified the three above functorial theories of characteristic classes in a functorial sense, similarly as F. Hirzebruch unified the three theories of cohomological characteristic classes: Chern classes, Todd classes and Thom-Hirzebruch $L$-classes. They proved the unique existence of a natural transformation 
$$T_{y,*}\colon K_0(var/-)\to H_{2*}(-,\Q)\otimes\Q[y],$$
 where $K_0(var/-)$ is the covariant relative Grothendieck functor of complex algebraic varieties, satisfying $T_{y,*}([X\xrightarrow{id_X}X])=T_y^*(TX)\cap[X]$ if $X$ is a smooth variety. The {\em Hirzebruch (homology) class} of $X$ is defined as $T_{y,*}(X):=T_{y,*}([X\xrightarrow{id_X}X])$. Furthermore, they proved the unique existence of three natural transformations $\epsilon,$ $mC_0$, and $sd$ defined from the relative  Grothendieck group of varieties as follows:  A transformation $\epsilon\colon K_0(var/-)\to F(-),$ such that  $\epsilon([X\xrightarrow{id_X}X])=1_X$ if $X$ is smooth, and $T_{-1,*}=c_*\circ \epsilon$; a transformation
 $mC_0\colon K_0(var/-)\to G_0(-),$ such that $mC_0([X\xrightarrow{id_X}X])=\mathcal O_X$  if $X$ is smooth, and $T_{0,*}=td_*\circ mC_0$; and a transformation 
 $$sd\colon K_0(var/-)\to \Omega_{\K}(-),$$ such that  $sd([X\xrightarrow{id_X}X])=[\mathbb K_X[\dim_\C X]]$  if $X$ is smooth, and $T_{1,*}=L_*\circ sd$. 
 
 This shows that  the transformation $T_{y,*}$ unifies the transformations $c_*$, $td_*$ and $L_*$ via the transformations $\epsilon$, $mC_0$ and $sd$, respectively. However, this does not necessarily mean that the Hirzebruch (homology) class $T_{y,*}(X)$ specializes to  $c_*^{SM}(X)$ for $y=-1$, $td_*^{BFM}(X)$ for $y=0$, and $L_*(X)$ for $y=1$,  for any singular variety $X$. Indeed, J.-P. Brasselet, J. Schürmann and S. Yokura proved that $T_{-1,*}(X)=c_*^{SM}(X)\otimes \mathbb Q$  for   $y=-1$;   $T_{0,*}(X)=td_*^{BFM}(X)$  if  $X$  has at most du Bois singularities for  $y=0$; and conjectured in \cite[Remark 0.1]{BSY:2010} the following statement for  $y=1$:
\begin{teo}\label{conj:BSY}
If $X$ is a compact complex algebraic variety that is a rational homology manifold, then
$$T_{1,*}(X)=L_*(X).$$
\end{teo}

It is important to note that for any compact variety $X$,  the class $T_{1,*}(X)$ can be different from the class $L_*(X)$ in general, see \cite[Remark 0.1, Example 3.1]{BSY:2010}. Furthermore, Theorem \ref{conj:BSY}  provides a complete answer about the relation between the Hirzebruch transformation $T_{y,*}\colon K_0(var/-)\to H_*(-)$ and Chern classes, Todd classes and $L$-classes for singular varities.

\smallskip
This conjecture was previously solved in the following special cases:  
\begin{itemize}
\item[$-$] In \cite{BSY:2010}, \cite{Ba:2019} for degree zero by Saito's Hodge index theorem for projective varieties.  
\item[$-$] In \cite{CMSS:2009} for certain hypersurfaces with isolated singularities.  
\item[$-$] In \cite{CMSS:2012} for $X = Y /G$, with $Y$ a projective $G$-manifold and $G$ finite group of agebraic autormorphisms. 
\item[$-$] In \cite{MS:2015} for projective simplicial toric varieties.  
\item[$-$] In \cite{Ba:2019} for normal projective complex $3$-folds at worst canonical singularities, trivial canonical divisor and  $\dim H^1(X,\mathcal O_X)> 0$.  
\end{itemize}

In 2020, the authors proved Theorem \ref{conj:BSY} for all projective varieties in the arxiv version of this article. Shortly thereafter, M. Saito together the authors gave a different proof of Theorem \ref{conj:BSY} for general compact varieties using the theory of mixed Hodge modules \cite{FPS:2023}.  
\smallskip

In \cite[Remark 5.4]{BSY:2010}, the following more general statement than Theorem \ref{conj:BSY} was also conjectured: Let $X$ be a compact complex algebraic variety, then
$$IT_{1,*}(X)=L_*(X),$$ where $IT_{y,*}(X):=MHT_{y,*}(IC_X^H[-\text{dim}_\C X])$ is the {\em Hodge-theoretic Hirzebruch class} of $X$, $MHT_{y,*}$  is a natural transformation from the Grothendieck functor of mixed Hodge modules, see \cite[Section 5]{BSY:2010}, and $IC_X^H$ is the intersection Hodge module on $X$. This conjecture can be regarded as a characteristic class generalization of Saito's Hodge index theorem that expresses the signature of a compact complex projective variety  in terms of the intersection Hodge numbers.  It is solved in the following special cases: 
\begin{itemize}
\item[$-$] In degree zero,  for the projective case by  Saito's Hodge index theorem \cite{BSY:2010},\cite{Ba:2019} and  for any compact variety \cite{FPS:2023}. 
\item[$-$] For all rational homology manifolds (Theorem \ref{conj:BSY}) \cite{FPS:2023}. 
\item[$-$] For all Schubert varieties in a Grassmanian  \cite{BSW:2023}.
\end{itemize}

\subsection{Main results}

The present note is the first, original proof of Theorem \ref{conj:BSY} for projective varieties. This is a geometric proof based only on cubical hyperresolutions, the Decomposition Theorem and classical Hodge theory. We remark that our proof does not use Hodge modules.
However, the significance of our work goes beyond the theorem itself; it lies in the methods we employ, each possessing intrinsic interest. We highlight two key results:

Firstly, we prove a new characterization of rational homology manifolds in terms of cubical hyperresolutions, this result corresponds to Theorem \ref{teo:equivalentRHM}. Its proof uses classical Hodge theory in the sense of de Cataldo and Migliorini's proof of the Decomposition Theorem, and it is the only point in our proof of Theorem \ref{conj:BSY} where the projective assumtion is needed as hyperplane sections are used. This result is announced in Section \ref{subsec:perverseexact} and proved in Section \ref{sec:exactperverse}.

Secondly, we prove that the cobordism class of certain perverse sheaves such that their stalks are polarizable pure Hodge structures is independent of the choice of the polarization, see Lemma \ref{lem:easydetermination}. This result is a much stronger statement than its counter part in \cite{FPS:2023}.  See Section \ref{subsec:polarization} and Section \ref{section:ComputationCobordismgroup} for more details and the proof of this statement.
\medskip

In this paper, using the two above main results, we prove the following identity of cobordism classes, which is a stronger statement than Theorem \ref{conj:BSY} also conjectured in \cite{BSY:2010}:

\begin{teo}
\label{main2}
If $X$ is a projective complex algebraic variety that is a rational homology manifold, then 
\begin{equation}\label{eq:main2}
sd_\R([X\xrightarrow{id_X} X])=[\mathbb R_X[\dim_{\mathbb C} X], S]\in \Omega_\R(X),
\end{equation}
where $S$ is the self-duality defined by multiplication without any sign.
\end{teo}

Here, $sd_\R$ denotes the natural transformation $sd$ taking $\K=\R$. This result has been also proved for general compact varieties in \cite{FPS:2023} using Hodge modules.  Notice that $sd_\R([X\xrightarrow{id_X} X])\neq [\mathbb R_X[\dim_{\mathbb C} X], S]$ in general \cite[Remark 0.1]{BSY:2010}. Additionally, it is important to notice that Theorem \ref{main2} is false for $\K=\Q$, see \cite{FPS:2023}.  Theorem \ref{conj:BSY} follows after applying the Cappell-Shaneson-Weinberger $L$-transformation $L_*$ to Equality (\ref{eq:main2}) and using the equality of transformations $T_{1,*}=L_*\circ sd_\R$.

\subsubsection{\textbf{\em Perverse exact cubical hyperresolutions}}\label{subsec:perverseexact} In~\cite{GNPP:1988} cubical hyperresolutions and their associated semi-simplicial resolutions were introduced. We follow the treatment of~\cite[Chapter 5]{PeSte:2008} and we refer to this source for the relevant definitions. \smallskip

Let $X$ be a complex algebraic variety of complex dimension $n$. Let  $\{Z_I\}_{I\subseteq [n]}$ be an $(n+1)$-cubical hyperresolution $X$ and let $\varepsilon\colon X_{\bullet}\to X$ be its associated semi-simplicial resolution. By \cite[Section 5.1.2]{PeSte:2008},  we have the following complex in $D_c^b(X)$: $$0\to R\varepsilon_*\R_{X_0} \to R\varepsilon_*\R_{X_1} \to \cdots \to R\varepsilon_*\R_{X_n} \to 0.$$ Applying the perverse cohomology functor ${^p\hH}^q(-)\colon D_c^b(X)\to Perv(X)$ to the above complex, for each $q\geq 0$, we obtain the following complex in the abelian category $Perv(X)$:
\begin{equation}\label{eq:complexperversenosup}
0\to {^p\hH}^q(R\varepsilon_*\R_{X_0}) \to {^p\hH}^q(R\varepsilon_*\R_{X_1}) \to \cdots \to {^p\hH}^q(R\varepsilon_*\R_{X_n}) \to 0.
\end{equation}

\begin{defi} \label{def:exactperverse} An $(n+1)$-cubical hyperresolution  $\{Z_I\}_{I\subseteq [n]}$ of $X$ with  associated semi-simplicial resolution $\varepsilon\colon X_{\bullet}\to X$ is {\em perverse exact} if the following properties hold:
\begin{enumerate}
\item[(i)] For every $q\neq n$, the complex (\ref{eq:complexperversenosup}) is exact;\smallskip
\item[(ii)] For $q=n$,  the perverse cohomology of the complex (\ref{eq:complexperversenosup}) is only concentrated in its 0-degree part and it is equal to $IC_X[-n]$.
\end{enumerate}
\end{defi}

\begin{re} \label{re:supcomplexperverse} Let $\{X=\Sigma_0,\Sigma_1,\dots,\Sigma_N\}$ be the collection of subvarieties in $X$ which are the support of the simple direct summands of the perverse sheaves ${^p\hH}^q(R\varepsilon_*\R_{X_k})$.  By the Decomposition Theorem (Theorem \ref{teo:DecTheo}), each perverse sheaf ${^p\hH}^q(R\varepsilon_*\R_{X_k})$ is decomposed as
\begin{equation}
\label{eq:decom001}
 {^p\hH}^q(R\varepsilon_*\R_{X_k})=\bigoplus_{j=0}^N{^p\hH}^q(R\varepsilon_*\R_{X_k})_{\Sigma_j},
\end{equation}
where ${^p\hH}^q(R\varepsilon_*\R_{X_k})_{\Sigma_j}$ denotes the direct sum of the simple summands in ${^p\hH}^q(R\varepsilon_*\R_{X_k})$ whose support is $\Sigma_j$. By the semi-simplicity of the perverse sheaves ${^p\hH}^q(R\varepsilon_*\R_{X_k})$, the complex (\ref{eq:complexperversenosup}) splits as the direct sum of the complexes
\begin{equation}
\label{eq:complexperverse}
0\to {^p\hH}^q(R\varepsilon_*\R_{X_0})_{\Sigma_j} \to {^p\hH}^q(R\varepsilon_*\R_{X_1})_{\Sigma_j} \to \cdots \to {^p\hH}^q(R\varepsilon_*\R_{X_n})_{\Sigma_j} \to 0,
\end{equation}
for $j=0,\dots,N$. 

Note that property  (i) of Definition  \ref{def:exactperverse} is equivalent to the exactness of complex (\ref{eq:complexperverse}) for all $p\neq n$ and all $\Sigma_j$. Furthermore, property (ii) of Definition  \ref{def:exactperverse}  is equivalent to the exactness of complex (\ref{eq:complexperverse}) for $q=n$ and all $\Sigma_j$ strictly contained in $X$, and  ${^p\hH}^{n}(R\varepsilon_*\R_{X_0})_{X}=IC_X[-n]$. 
\end{re}

We call a projective of degree 1 cubical hyperresolution of $X$ to a cubical hyperresolution which has only one component of dimension $n$ and such that this component is a resolution of singularities of $X$, see Definition \ref{def:degcubicalhyp} for more details. We will use these type of cubical hyperresolutions to prove Theorem \ref{main2}. Notice that this is the type of cubical hyperresolutions constructed in \cite{GNPP:1988}, \cite{PeSte:2008}, see Remark \ref{re:PeStecub}.

\medskip
We prove the following characterization theorem of rational homology manifolds in terms of perverse exact cubical hyperresolutions.

\begin{teo}\label{teo:equivalentRHM} Let $X$ be a projective variety of complex dimension $n$. The following statements are equivalent:
\begin{enumerate}
\item[(i)] $X$ is a rational homology manifold;
\item[(ii)] There exists a projective of degree 1 cubical hyperresolution of $X$ that is perverse exact;
\item[(iii)] All projective of degree 1 cubical hyperresolutions of $X$  are perverse exact.
\end{enumerate}
\end{teo}

One significant application of this result is in proving the vanishing of the cobordism difference class  $sd_\R([X\xrightarrow{id_X}])-[IC_X]\in \Omega_\R(X)$ when $X$ is a rational homology manifold, as detailed in Section \ref{section:ComputationCobordismgroup}. The proof of Theorem \ref{teo:equivalentRHM} is given in Section \ref{sec:exactperverse}.

\subsubsection{\textbf{\em Independence of polarizations in the cobordism group}} \label{subsec:polarization} Another essential intermediate step, needed both in \cite{FPS:2023} and here, is the fact that the isomorphism class of a pair $(\mathcal F, S)$ given by some constructible complex of sheaves $\mathcal F$ together with a perfect pairing $S\colon \mathcal F\otimes\mathcal F\to \mathbb D_X$  is independent of the choice of pairing $S$.  This holds true as long as the pairing is a polarization of Hodge modules or a polarization of Hodge structures at the stalks of local systems, respectively.  In particular, this assertion holds in the cobordism group of self-dual $\R$-complexes of sheaves. In this paper this result corresponds to Lemma~\ref{lem:easydetermination}, and Proposition 1 in \cite{FPS:2023}. 

We would like to emphazise  that Lemma~\ref{lem:easydetermination} is much stronger than \cite[Proposition 1]{FPS:2023}  as it does not require the semi-simplicity property on the local systems involved; see also \cite[Remark 2.1]{FPS:2023}. Consequently, the proof of Lemma~\ref{lem:easydetermination} diverges from that of~\cite[Proposition 1]{FPS:2023}.  To prove Lemma~\ref{lem:easydetermination} we use Proposition \ref{prop:rep} which shows that given a smooth 1-parameter family of non-degenerate symmetric or anti-symmetric parings $Q_s$ of a $\R$-local system $\mathcal L$, for $s$ in an interval $I$, one can find a smooth 1-parameter family of automorphisms $T(s)$ of $\mathcal L$, locally on $I$, such that the pairings $Q_s(T(s)(-),T(s)(-))$ are independent of the parameter $s$.

The proof of Proposition \ref{prop:rep} is based on a representation theory argument and integral curves of time dependent vector fields. Both Lemma~\ref{lem:easydetermination} and Proposition \ref{prop:rep} may be of some utility elsewhere. Choosing the field $\R$ is not a coincidence since one needs to take square roots in the proof of \cite[Proposition 1]{FPS:2023} and in Lemma \ref{lem:easydetermination}.



\subsection{Background} \label{sec:background} 
Let us precise what we mean by ``a geometric proof based on classical Hodge theory''.  Our proof uses the Decomposition Theorem for $\R$-coefficients in the following form: 

\begin{teo}[\cite{BBD:1982},~\cite{Saito:1989a},~\cite{DeCM:2005}]\label{teo:DecTheo}
Let $\varepsilon:Z\to X$ be proper morphism of complex algebraic varieties with $Z$ smooth of dimension $d$. Then,
\begin{equation}
\label{eq:d00}
R\varepsilon_*\R_Z[d]\cong \bigoplus_{i=-M}^M {^p\hH}^i(R\varepsilon_*\R_{Z}[d])[-i],
\end{equation}
where $M$ is a positive integer, and ${^p\hH}^i(-)$ denotes the $i$-th cohomology functor for the perverse $t$-structure introduced in~\cite{BBD:1982}. Moreover, the perverse sheaves ${^p\hH}^i(R\varepsilon_*\R_{Z}[d])$ are semi-simple. 
\end{teo} 

The Decomposition Theorem was proved originally in~\cite{BBD:1982} for $\C$-coefficients. For $\R$-coefficients, the equation~(\ref{eq:d00}), and even its analogue for $\Q$-coefficients follows from~\cite{De:1968}; then the semi-simplicity of the perverse sheaves ${^p\hH}^q(R\varepsilon_*\R_{Z})$ can be obtained by an argument involving the exactness of the scalar extension functor from the abelian category of $\R$-perverse sheaves to the abelian category of $\C$-perverse sheaves. We thank M. Saito for this remark. The Decomposition Theorem in the form we need was re-proved later by  M. Saito, using the theory of Hodge modules (see \cite{Saito:1989},~\cite{Saito:1989a},~\cite{Saito:1990}), and also by M. A. de Cataldo and L. Migliorini in~\cite{DeCM:2005}. M. A. de Cataldo and L. Migliorini's proof is geometric and rests in classical Hodge theory in the sense that it only uses the formalism of perverse sheaves and Hodge theory as developed in~\cite{Deligne:1971},~\cite{Deligne:1974}. 


Given the Decomposition Theorem as stated above, the theory of cubical hyperresolutions and generalities in perverse sheaves, our proof  of Theorem \ref{main2} only needs fairly elementary Hodge theory computations as developed in~\cite{Deligne:1971},~\cite{Deligne:1974}.


\subsection*{Conventions}
(i) In this paper, we denote simply by $[X']$ the class $[X'\to X]\in K_0(var/X)$ if there is no risk of confusion.

(ii) We use Deligne's indexing convention for complexes of sheaves and self-dualities. If $X$ is smooth variety of complex dimension $n$, then the dualizing complex of $X$ is $\mathbb D_X=\R_X[2n]$ and the intersection cohomology sheaf complex is $IC_X=\R_X[n]$. Usual real numbers multiplication (without any sign) defines a  pairing $\sigma_X:\R_X[n]\otimes \R_X[n]\to \R_X[2n]$, which induces a self-duality in $IC_X=\R_X[n]$. Note that this convention differs from the one in \cite{GM2:1983} and \cite{CS:1991}.

(iii) In \cite[Section 1.1]{FPS:2023} a new definition of cobordism group of self-dual complexes of sheaves has been introduced. This definition uses perfect pairings rather than self-dualities (see Section \ref{section:ComputationCobordismgroup}). In alignment with \cite{FPS:2023}, we adopt this definition of cobordism group in this paper.  The elements in $\Omega_{\K}(X)$ are equivalence classes $[ \mathcal F,S]$ of pairs $(\mathcal F,S)$ where $\mathcal F$ is a self-dual complex of sheaves on $X$ and $S$ is a self-duality on $\mathcal F$. Sometimes we will does not specify the self-duality $S$, and we will write only $[\mathcal F]$ instead of $[\mathcal F,S]$.

(iv) In this paper, we follow the definition of polarization as it appears in~\cite[Definition 2.9]{PeSte:2008}.

\subsubsection*{Acknowledgements}The first author thanks M. Banagl for a beautiful course that raised his interest in the conjecture. We also thank to L. Maxim, V. Mu\~noz J. Sch\"urmann and S. Yokura for useful comments and positive  criticism that improved the paper a lot. M. Saito generous help was very important, at several points.


\section{An identity in the Grothendieck group of varieties}\label{sec:groth}

In this section we give an identity in the relative Grothendieck group of varieties $K_0(var/X)$ that expresses  the class of the identity map  $[X\xrightarrow{id_X}X]\in K_0(var/X)$ as an alternate sum of classes from smooth varieties to $X$  coming from a cubical hyperresolution. This expression allows to compute $sd_\R([X\xrightarrow{id_X}X])\in \Omega_\R(X)$ such that the cobordism class of the intersection cohomology sheaf complex $IC_X$ appears in the expression together with other terms, as we show in Section \ref{section:ComputationCobordismgroup}. 




\medskip
Let $X$ be a complex algebraic variety of complex dimension $n$. Set $[n]:=\{0,\dots,n\}$.

\begin{defi}  \label{def:degcubicalhyp} An $(n+1)$-cubical hyperresolution  $\{Z_I\}_{I\subseteq [n]}$ of $X$ is {\em projective of degree 1} if its associated  semi-simplicial resolution $\varepsilon\colon X_{\bullet}\to X$ satisfies the following properties:
\begin{itemize}
\item[(i)] Each $X_k$ is a disjoint union of smooth varieties $X_{k,i}$, such that $\dim X_{k,i}\leq n$; 
\item[(ii)] The term $X_0$ is of the form $X_0=\tilde{X}\bigsqcup (\bigsqcup_{i=1}^n X_{0,i})$ where its only compoment of dimension $n$ is the smooth variety $\tilde{X}$ and $\varepsilon|_{\tilde{X}}:\tilde{X}\to X$ is a resolution of singularities of $X$.
\item[(iii)] All the morphisms $X_k\to X$ are projective.
\end{itemize}
\end{defi}

\begin{re}\label{re:PeStecub}
The cubical hyperresolutions constructed in \cite[Section 5.2, Theorem 5.26]{PeSte:2008} are projective of degree 1. We also note that all morphisms involved in the cubical hyperresolution are projective.
\end{re}

Let $K_0(var/X)$ be the relative Grothendieck group of varieties over $X$. There is a cubical hyperresolution with associated semi-simplicial resolution $\varepsilon\colon X_{\bullet}\to X$ that is projective of degree 1, such that the following identity holds in $K_0(var/X)$:
\begin{equation}\label{Eq1}
[X]=[\tilde X]+[\bigsqcup_{i=1}^n X_{0,i}]+\sum_{k=1}^n (-1)^k[X_k].
\end{equation}

The identity (\ref{Eq1}) follows by considering the projective of degree 1 cubical hyperresolution constructed in the proof of \cite[Theorem~5.26]{PeSte:2008}, see Remark \ref{re:PeStecub}, observing that the corresponding discriminant square \cite[(V-11)]{PeSte:2008} at each step of the construction is cartesian, and applying the additivity relation of $K_0(var/X)$ to it.

\begin{lem}
\label{rem:hypsection}
Let $\{Z_I\}_{I\subseteq [n]}$ be a projective of degree 1 cubical hyperresolution of a projective variety $X\subset \PP^N_{\C}$. Let $H\in (\PP^N_{\C})^*$ be a generic hyperplane. Then, the cubical scheme formed taking the fibre product of $\{Z_I\}_{I\subseteq [n]}$ by $X\cap H$ is a projective of degree 1 cubical hyperresolution of $X\cap H$. Furthermore, a variation $H_t$ of the generic hyperplane where $H_t\in U\subset (\PP^N_{\C})^*$ is a small neighborhoud of the point $H_0$ in the dual projective space $(\PP^N_{\C})^*$, yields a topologically trivial family of cubical hyperresolutions. Iterating we obtain the same statement for generic linear sections of arbitary codimension.
\end{lem}
\proof
The proof of the first assertion is an inspection on the construction of the cubical hyperresolution in~\cite[Section 5.2]{PeSte:2008}, combined with the fact that if $|L|\subset (\PP^M_{\C})^*$ is a linear system without base points in a smooth projective manifold $Z\subset\PP^M_\C$, then a generic hyperplane section $Z\cap H$ in $|L|$ is smooth. 
\endproof  



\section{Proof of Theorem \ref{teo:equivalentRHM}}\label{sec:exactperverse}

In Theorem \ref{teo:equivalentRHM}, we establish a new characterization of rational homology manifolds using cubical hyperresolutions.  In Section \ref{subsec:spectralsec}, we construct two spectral sequences of perverse sheaves associated to a cubical hyperresolution. We will use them to prove Theorem \ref{teo:equivalentRHM} in Section \ref{subsec:proofequivalent}. 

\subsection{Rational homology manifolds}\label{RHM}
\begin{defi} A complex algebraic variety $X$ of complex dimension $n$ is a {\em rational homology manifold} if, for every $x\in X$, we have
$$H_i(X,X\setminus x;\Q)=\left\{ \begin{array}{lc} \Q, & i =2n, \\ \\ 0, &  i\neq 2n. \end{array} \right.$$
\end{defi}

\begin{propo}[\cite{BM:1983}, \cite{M:2009}]  \label{propo:RHMchar} A complex algebraic variety $X$ of complex dimension $n$ is a rational homology manifold if, and only if, the natural morphism $\mathbb R_X[n]\to IC_X$ is an isomorphism in $D_c^b(X)$. 
\end{propo}

To prove Theorem \ref{teo:equivalentRHM}, we require the following lemma regarding hyperplane sections.
\begin{lem}
\label{lem:qhmsection} 
Any generic hyperplane section of a projective variety $X$ that is a rational homology manifold is also a rational homology manifold.
\end{lem}
\proof
Choose a Whitney stratification of $X$. Note that a  generic hyperplane $H$  does not meet the $0$-dimensional stratum, and so for any point $x\in X\cap H$,  there exists a neighborhood $U$ of $x$ in $X$ such that $U=(U\cap H)\times D$, where $D$ is a disk. The proof follows now from easy homological considerations.
\endproof

\subsection{Two spectral sequences of perverse sheaves} \label{subsec:spectralsec}
Let $X$ be a complex algebraic variety of dimension $n$ and let $\varepsilon\colon X_\bullet \to X$ be the associated semi-simplicial resolution of an $(n+1)$-cubical hyperresolution of $X$.

\begin{lem} \label{lem:spectralsecI}
There exists a spectral sequence of perverse sheaves associated to $\varepsilon\colon X_\bullet \to X$ such that the complexes of its $E_1$-page  coincide with the complexes (\ref{eq:complexperversenosup}).
\end{lem}

\begin{proof}
Denote by $\mathcal C^\bullet_{X_k}$ the canonical Godement resolution of the constant sheaf $\R_{X_k}$ in $X_k$ for each $k$. We obtain a double complex of sheaves $I^{\bullet,\bullet}$ in $X$ associated to the semi-simplicial resolution $\varepsilon\colon X_{\bullet}\to X$,  such that each column $I^{k,\bullet}$ is equal to $\varepsilon_*\mathcal C^\bullet_{X_k}$, and hence computes $R\varepsilon_*\R_{X_k}$. It is important to notice that the horizontal differentials in the double complex $I^{\bullet,\bullet}$ give rise to morphisms of complexes $\varepsilon_*\mathcal C^\bullet_{X_k}\to \varepsilon_*\mathcal C^\bullet_{X_{k+1}}$ which are induced from the alternating sum of the pullbacks by the $(k+1)$-morphisms $X_{k+1}\to X_{k}$ appearing in the semi-simplicial scheme $X_\bullet$ (see~\cite[Section 5.1.2]{PeSte:2008} for more details).

The single complex $s(I^{\bullet,\bullet})$ is decreasingly filtered by the subcomplexes $F^ps(I^{\bullet,\bullet})$, where $F^ps(I^{\bullet,\bullet})$ is the single complex of the double sub-complex of $I^{\bullet,\bullet}$ formed by the direct sum of $I^{a,b}$ for $a\geq p$. 

Define the Cartan-Eilenberg systems of perverse sheaves as follows (see~\cite[XV.7]{CaEi:1956}, which makes sense for the abelian category of perverse sheaves):
\begin{equation}\label{sum}
H[I^{\bullet,\bullet}](p,q):=\sum_{i\in\Z}{^p\hH}^i(F^ps(I^{\bullet,\bullet})/ F^qs(I^{\bullet,\bullet})),
\end{equation}
for $p\geq q$. The morphisms $H[I^{\bullet,\bullet}](p,q)\to H[I^{\bullet,\bullet}](p',q')$ are induced from the natural morphism of complexes. The connecting morphisms coincide with the connecting morphism for the exact sequence of complexes 
$$0 \to F^qs(I^{\bullet,\bullet})/ F^rs(I^{\bullet,\bullet})\to F^ps(I^{\bullet,\bullet})/ F^rs(I^{\bullet,\bullet})\to F^ps(I^{\bullet,\bullet})/ F^qs(I^{\bullet,\bullet})\to 0$$ The defining sums above are finite, and it is straightforward to check, using that $X_k=\emptyset$ for $k>n$, that conditions (SP.1)-(SP.5) of \cite[XV.7]{CaEi:1956} are satisfied (for (SP.5) it is convenient to use that $X_k$ is non empty for finitely many $k$). The perverse sheaves $H[I^{\bullet,\bullet}](p,q)$ are graded by the defining sum (\ref{sum}).

Therefore, by~\cite[XV.7]{CaEi:1956}, we obtain an spectral sequences of graded objects in $Perv(X)$. We obtain the following terms: 
\begin{equation}
\label{eq:E1}
 E(I)_1^{p,q}\cong {^p\hH}^q(R\varepsilon_*\R_{X_p}).
\end{equation}

Then, the complexes appearing in the $E_1$-page of the spectral sequence associated with the complex $s(I^{\bullet,\bullet})$ coincide with the complexes~(\ref{eq:complexperversenosup}).   
\end{proof}

Furthermore, by the semi-simplicity of the perverse sheaves predicted by the Decomposition Theorem (Theorem \ref{teo:DecTheo}), the complexes appearing in the $E_1$-page of the above spectral sequence of perverse sheaves split in a direct sum of complexes of perverse sheaves with strict support $\Sigma_j$.  These complexes coincide with the complexes~(\ref{eq:complexperverse}).

\begin{lem}  \label{lem:spectralsecK} There exists a spectral sequence of perverse sheaves defined by a single complex $K^\bullet$ associated to a cubical hyperresolution that is projective of degree 1, satisfying that the complexes of its $E_1$-page coincide with the complexes (\ref{eq:complexperverse}) for all supports strictly contained in $X$. Furthermore,  $X$ is a rational homology manifold if, and only if, the complex $K^\bullet$ is acyclic.
\end{lem}
\begin{proof}
Let  $I^{\bullet,\bullet}$ be the double complex appearing in the proof of Lemma \ref{lem:spectralsecI} associated to the semi-simplicial resolution $\varepsilon\colon X_\bullet \to X$.  Consider the sequence of morphisms of double complexes 
$$I^{\bullet,\bullet}\stackrel{\beta_1}{\longrightarrow}I^{0,\bullet}\stackrel{\beta_2}{\longrightarrow}R(\varepsilon|_{\tilde{X}})_*\R_{\tilde{X}}\stackrel{\beta_3}{\longrightarrow}IC_X[-n],$$
where $I^{0,\bullet}$ is the double complex that has non-zero objects only at the $0$-th column, $\beta_1$ is the natural morphism of double complexes, $\beta_2$ is the composition of the natural projection from $I^{0,\bullet}$ to $R(\varepsilon|_{\tilde{X}})_*\R_{\tilde{X}}$ given by the decomposition of $X_0$ in connected components, and $\beta_3$ is described as follows: by Theorem \ref{teo:DecTheo} there is a non-canonical direct sum decomposition
$$\Phi:R(\varepsilon|_{\tilde{X}})_*\R_{\tilde{X}}\to IC_X[-n]\oplus L,$$
where $L$ is a direct sum of shifted simple perverse sheaves with support strictly contained in $X$. We define $\beta_3:=\rho\circ\Phi$, where $\rho$ is the canonical projection $IC_X[-n]\oplus L\to IC_X[-n]$.

Since $\varepsilon \colon X_\bullet\to X$ is of cohomological descent \cite[Definition 5.6]{PeSte:2008}, there is a quasi-isomorphism  $\R_X\to s(I^{\bullet,\bullet})$. We define the single complex morphism 
\begin{equation}
\label{eq:eta}
\eta:s(I^{\bullet,\bullet})\cong\R_X\to IC_X[-n]
\end{equation}
as the composition  $\beta_3\circ\beta_2\circ\beta_1$. Even if $\beta_3$ is not unique, by the uniqueness of Proposition in~\cite[Section 5.1]{GM2:1983}, the morphism $\eta$ is, up to multiplication with a non-zero real number, the canonical morphism connecting cohomology with intersection cohomology complexes.

Notice that we can form a double complex $K^{\bullet,\bullet}$ whose columns are $K^{0,\bullet}=cone(\beta_3\circ\beta_2)[-1]$ and $K^{p,\bullet}=I^{p,\bullet}$ for $p>0$ (in this formula $\beta_3\circ\beta_2$ denotes the simple complex morphism induced at the $0$-th column).  Furthermore,  there is a quasi-isomorphism of single complexes 
\begin{equation}\label{eq:qisconeK}
s(K^{\bullet,\bullet})\cong cone(\eta)[-1].
\end{equation}
We notice that the double complex $K^{\bullet,\bullet}$ depends on the choice of $\beta_3$. However, this non-uniqueness does not affect our proof.  

By Proposition \ref{propo:RHMchar}, $X$ is a rational homology manifold if, and only if, the morphism $\eta$  is a quasi-isomorphism i.e. the $cone(\eta)[-1]$ is acyclic. By the quasi-isomorphism (\ref{eq:qisconeK}), $X$ is a rational homology manifold if, and only if, the single complex $s(K^{\bullet,\bullet})$ is acyclic.
\medskip

Similar to the construction of the spectral sequence of perverse sheaves for the single complex $s(I^{\bullet,\bullet})$ in Lemma \ref{lem:spectralsecI} using Cartan-Eilenberg systems \cite[XV.7]{CaEi:1956}, we obtain a spectral sequence of perverse shaves for the single complex $s(K^{\bullet,\bullet})$. Moreover, there is an obvious morphism of Cartan-Eilenberg systems $H[I^{\bullet,\bullet}](p,q)\to H[K^{\bullet,\bullet}](p,q)$ and a morphism between both spectral sequences.  Since $K^{p,q}\to I^{p,q}$ is an isomorphism, for $p>0$, we have an isomorphism 
\begin{equation}
\label{eq:E1p}
 E(K)_1^{p,q}\cong E(I)_1^{p,q}\cong {^p\hH}^q(R\varepsilon_*\R_{X_p}).
\end{equation}
Considering the decomposition~(\ref{eq:decom001}), and by definition of $K^{\bullet,\bullet}$, for $p=0$ we have
\begin{equation}
\label{eq:E10}
 E(K)_1^{0,q}:=\bigoplus_{j\neq 0}{^p\hH}^q(R\varepsilon_*\R_{X_0})_{\Sigma_j};
\end{equation}
that is, all the summands except $IC_X[-n]$ if $q=n$. 

By the semi-simplicity of the perverse sheaves predicted by Theorem \ref{teo:DecTheo}, the complexes appearing in the $E_1$-page of the spectral sequence associated with the complex $s(K^{\bullet,\bullet})$ splits in a direct sum of complexes of perverse sheaves with support $\Sigma_j$  which coincides with the complex (\ref{eq:complexperverse}) with the exception of the term ${^p\hH}^n(R\varepsilon_*\R_{X_0})_X$.  
\end{proof}

\subsection{Proof of Theorem \ref{teo:equivalentRHM}} \label{subsec:proofequivalent}


We prove the theorem in two steps: (ii) $\longrightarrow$  (i) and (i) $\longrightarrow$  (iii).  The hard part of the proof is (i) $\longrightarrow$  (iii), which we achieve by proving the degeneration at the $E_2$-page of the spectral sequence of perverse sheaves constructed in Lemma \ref{lem:spectralsecK}.  Our proof only uses classical Hodge theory (see Section \ref{sec:background}), but it applies only when $X$ is projective due to our use of hyperplane sections. However, the part (ii) $\longrightarrow$  (i) holds for any complex algebraic variety $X$.  We notice that  Theorem \ref{teo:equivalentRHM} is the only point in the proof of Theorem \ref{main2} that the projectivity assumtion on $X$ is needed. 

(ii) $\longrightarrow$  (i). 
Suppose that exists a projective of degree 1 cubical hyperresolution that is perverse exact and let $\varepsilon\colon X_\bullet\to X$ its associated semi-simplicial resolution. By the exactness of the complexes and by the proof of Lemma \ref{lem:spectralsecK}, the spectral sequence of perverse sheaves defined by the single complex $s(K^{\bullet, \bullet})$ degenerates at the $E_2$-page, and the terms $E_2^{p,q}$ vanish for all $p,q$. Hence, the single complex $s(K^{\bullet, \bullet})$ is acyclic, so $X$ is a rational homology manifold.
\medskip

(i) $\longrightarrow$ (iii). Assume that $X$ a rational homology manifold and choose any projective of degree 1 cubical hyperresolution. We prove that it is perverse exact. Notice that proving the degeneration at the $E_2$-page of the spectral sequence of perverse sheaves constructed in the proof of Lemma \ref{lem:spectralsecK} associated to our cubical hyperresolution is enough to prove the exactness of the complexes~(\ref{eq:complexperverse}) for $\Sigma_j$ strictly contained in $X$, which implies properties (i) and (ii) of Definition \ref{def:exactperverse}.

The proof of the degeneration at $E_2$ is by double induction on $\dim X$ and $\text{codim } \Sigma_j$. Suppose that the statement holds for $\dim X<n$ and for $\text{codim }\Sigma_j<d$ when $\dim X=n$. Assuming $\dim X=n$, we prove the exactness simultaneously for all supports $\Sigma_j$ of codimension $d$ in $X$. 
\medskip

\textsc{Case $d<n$}:  For any $\Sigma_j$ there exists a dense open subset $U_j$ over which all the perverse sheaves $(E_{1}(K)^{p,q})_{\Sigma_j}$ are local systems. In order to prove the exactness of (\ref{eq:complexperverse}) it is enough to prove the exactness at the stalk
\begin{equation}
\label{eq:shortexactstalk}
0\to((E(K)_{1}^{0,q})_{\Sigma_j})_z\to ((E(K)_{1}^{1,q})_{\Sigma_j})_z\to \cdots \to ((E(K)_{1}^{n,q})_{\Sigma_j})_z\to 0,
\end{equation}
of the complex~(\ref{eq:complexperverse}) at a point $z$ in each connected component of each of the open subsets $U_j$. Let $H$ be a generic linear section of dimension $d$ such that the intersection $\Sigma_j\cap H$ is a finite set of points contained in $U_j$ for every $d$-codimensional component $\Sigma_j$. Then, $X\cap H$ is a projective rational homology manifold of dimension $n-\dim \Sigma_j$ by Lemma~\ref{lem:qhmsection}, and by Lemma~\ref{rem:hypsection} the pullback to $X\cap H$ of the semi-simplicial resolution of $\varepsilon\colon X_{\bullet}\to X$ gives a semi-simplicial resolution of $(X|_H)_\bullet\to X\cap H$. Construct the perverse spectral sequence of the proof of Lemma \ref{lem:spectralsecK} for the   semi-simplicial resolution $(X|_H)_\bullet\to X\cap H$, and split it as a direct sum of spectral sequences of perverse sheaves with common support as above. 

For any $z\in \Sigma_j\cap H$ the point $z$ is a support for the $E_1$-page of the perverse spectral sequence associated with the hyperresolution $(X|_H)_\bullet\to X\cap H$ and the complex~(\ref{eq:shortexactstalk}) is the analog of the complex~(\ref{eq:complexperverse}) for the support $z$. 
This follows because $(X|_H)_p$ is the fibre product $X_p\times_{X} (X\cap H)$, and then, by the topological triviality statement of Lemma~\ref{rem:hypsection}, we have $R\varepsilon_*\R_{(X|_H)_p}=\iota_{X\cap H}^*R\varepsilon_*\R_{X_p}$, where $\iota_{X\cap H}$ denotes the inclusion of $X\cap H$ into $X$.

Since $\dim X\cap H =n-\dim\Sigma_j<n$, by induction hypothesis the result is true for $X\cap H$ and the semi-simplicial resolution $(X|_H)_\bullet\to  X\cap H$, we have the exatness of the sequence~(\ref{eq:shortexactstalk}). 
\medskip

\textsc{Case $d=n$}: Applying the cohomology functor $H^*(X,-)$ to the quotients $F^ps(I^{\bullet,\bullet})/ F^qs(I^{\bullet,\bullet})$ and $F^ps(K^{\bullet,\bullet})/ F^qs(K^{\bullet,\bullet})$ appearing in the proof of Lemma \ref{lem:spectralsecI} and Lemma \ref{lem:spectralsecK}, we obtain two Cartan-Eilenberg systems and a morphism between them, in a similar way as above for the construction of spectral sequences of perverse sheaves. They induce spectral sequences of real vector spaces, whose terms  are denoted by $'E(K)_r^{p,b}$ and $'E(I)_r^{p,b}$. The morphism between the Cartan-Eilenberg systems induces homomorphisms $'E(K)_r^{p,b}\to {'E(I)}_r^{p,b}$ which are compatible with the differentials. The $E_1$ terms are the following:
\begin{equation}
\label{eq:E1global}
'E(I)_1^{p,b}\cong H^b(X,R\varepsilon_*\R_{X_p}).
\end{equation}
For $p>0$ we have an isomorphism
\begin{equation}
\label{eq:E1pglobal}
'E(K)_1^{p,b}\cong {'E(I)}_1^{p,b}\cong H^b(X,R\varepsilon_*\R_{X_p}),
\end{equation}
and, for $p=0$ we have
\begin{equation}
\label{eq:E10global}
'E(K)_1^{0,b}:=\Ker (H^b(X,R\varepsilon_*\R_{X_0})\stackrel{\beta_3\circ\beta_2}{\longrightarrow} H^b(X,IC_X[-n])).              
\end{equation}

The spectral sequence $'E(I)$ coincides with the spectral sequence induced by the filtration by columns of the double complex $\Gamma(X,I^{\bullet,\bullet})$. By \cite[Theorem 3.18,  Theorem 5.33]{PeSte:2008}, $'E(I)$ lifts to a spectral sequence of real mixed Hodge structures, degenerates at $E_2$, converges to the mixed Hodge structure $H^*(X;\R)$, and for every $k$ we have that $W_{k-r}H^k(X;\R)\cong \bigoplus_{p\geq r}\mbox{} 'E(I)_2^{p,b}$. Since $X$ is a rational homology manifold, $H^k(X;\R)$ is a pure Hodge structure of weight $k$ \cite{Ze1:1983}, \cite{Ze2:1983}, and so $'E(I)_2^{p,b}=0$ for $p\geq 1$.

By the isomorphism~(\ref{eq:E1pglobal}) we deduce that $'E(K)_2^{p,b}\cong {'E(I)}_2^{p,b}=0$ for $p\geq 2$. Therefore $'E(K)_r^{p,b}=0$ for any $r$ and $p\geq 2$, and the spectral sequence $'E(K)$ degenerates at the $E_2$-page. 



By the degeneration at the second page, since the complex $s(K^{\bullet,\bullet})$ is acyclic, we have the vanishing $'E(K)_2^{p,b}=0$ for every $p,b$, and therefore for every $b$ we have the exact sequence
\begin{equation}
\label{eq:shortexactglobal}
0\to{'E(K)}_{1}^{0,b}\to {'E(K)}_{1}^{1,b}\to \cdots \to {'E(K)}_{1}^{n,b}\to 0.
\end{equation}
We denote by $d'_1$ the differential of this complex.  

By Theorem \ref{teo:DecTheo} and the $E_1$ term description~(\ref{eq:E1pglobal}) and~(\ref{eq:E10global}), we have
the splitting 
\begin{equation}
\label{eq:splitting0}
'E(K)_1^{p,b}=\bigoplus_{1\leq j\leq N}\bigoplus_{q\geq 0} H^b(X,(E(K)_1^{p,q})_{\Sigma_j}[-q])
\end{equation}
in the category of real vector spaces.

Denote by 
$$d'_1(p,b,j_1,j_2,q_1,q_2)\colon H^b(X,(E(K)_1^{p,q_1})_{\Sigma_{j_1}}[-q_1])\to H^b(X,(E(K)_1^{p+1,q_2})_{\Sigma_{j_2}}[-q_2])$$
the composition of $d'_1$ with the inclusion of the source in $'E(K)_1^{p,b}$ and the projection from $'E(K)_1^{p+1,b}$ to the target. For $j_1=j_2$, $q_1=q_2$ the morphism $d'_1(p,b,j_1,j_1,q_1,q_1)$ is obtained applying the functor $H^b(X,-)$ to the differential $d_1:(E(K)_1^{p,q_1})_{\Sigma_{j_1}}\to (E(K)_1^{p+1,q_1})_{\Sigma_{j_1}}$ appearing in the complex (\ref{eq:complexperverse}) for $q=q_1$. For $j_1\neq j_2$, $q_1=q_2$ the morphism $d'_1(p,b,j_1,j_2,q_1,q_1)$ vanishes since $(E(K)_1^{p,q_1})_{\Sigma_{j_1}}$ and $(E(K)_1^{p,q_1})_{\Sigma_{j_2}}$ are semi-simple perverse sheaves of disjoint support.

By induction, if $\Sigma_j\neq X$ is not of dimension $0$,  the sequence of semi-simple perverse sheaves~(\ref{eq:complexperverse}) is exact. Semi-simplicity and exactness imply that the sequence~(\ref{eq:complexperverse}) is isomorphic to the direct sum of exact sequences of perverse sheaves of the form $0\to P[-l]\to P[-(l+1)]\to 0$, where $P$ is simple and perverse, the map is the identity and $0\leq l\leq n-1$. 

Pick some $\Sigma_j\neq X$, and some simple perverse sheaf such that $0\to P[-l]\to P[-(l+1)]\to 0$ is a direct summand of~(\ref{eq:complexperverse}). Pick any $b$ such that $H^b(X,P[-q])\neq 0$. Then the sequence~(\ref{eq:shortexactglobal}) is isomorphic to one of the form
\begin{equation}
\label{eq:complejoauxiliar}
...\to A^{l-1}\to A^l\oplus H^b(X,P[-q])\stackrel{d'_1}{\longrightarrow} A^{l+1}\oplus H^b(X,P[-q])\to A^{l+2}\to \cdots
\end{equation}
where the differential $d'_1$ restricted and projected to $H^b(X,P[-q])$ is the identity. Then the sequence 
$$\cdots \to A^{l-1}\to A^l \to (A^{l+1}\oplus H^b(X,P[-q]))/d'_1(H^b(X,P[-q]))\to A^{l+2}\to\cdots,$$
is also exact, and identifying $A^{l+1}\cong (A^{l+1}\oplus H^b(X,P[-q]))/d_1(H^b(X,P[-q])$ is isomorphic to 
$$\cdots\to A^{l-1}\to A^l\to A^{l+1}\to A^{l+2}\to \cdots$$
with differential induced from $d'_1$. 

Set $d_j:=\dim \Sigma_j$. We proceed in the same way with all the direct summands $0\to P[-l]\to P[-(l+1)]\to 0$ appearing in all the exact sequences (\ref{eq:complexperverse}) for any $\Sigma_j\neq X$ not of dimension $0$ and any $q$. Taking into account the splitting~(\ref{eq:splitting0}), from the exact sequence~(\ref{eq:shortexactglobal}), we obtain an exact sequence of the form
$$\cdots\to \bigoplus_{d_j=0}\bigoplus_{q\geq 0} H^b(X,(E(K)_1^{p,q})_{\Sigma_j}[-q])\to \bigoplus_{d_j=0}\bigoplus_{q\geq 0} H^b(X,(E(K)_1^{p+1,q})_{\Sigma_j}[-q])\to \cdots $$
The group $H^b(X,(E(K)_1^{p,q})_{\Sigma_j})[-q]$ vanishes unless $q=b$, because $(E(K)_1^{p,q})_{\Sigma_j}$ has $0$-dimensional support, so the sequence becomes 
$$\cdots \to \bigoplus_{d_j=0} H^q(X,(E(K)_1^{p,q})_{\Sigma_j}[-q])\to \bigoplus_{d_j=0} H^q(X,(E(K)_1^{p+1,q})_{\Sigma_j}[-q])\to \cdots $$
Using that the morphism $d'_1(p,b,j_1,j_2,q_1,q_1)$ vanishes if $j_1\neq j_2$, we obtain that this sequence is the direct sum over the set of $0$-dimensional supports (for fixed $q$) of the sequences (\ref{eq:complexperverse}).

\section{A computation in the cobordism group of self-dual complexes}\label{section:ComputationCobordismgroup}

In this section we study self-dual $\R$-complexes of sheaves in the cobordism group $\Omega_\R(X)$ on a complex algebraic variety $X$. The cobordism group $\Omega_\R(X)$ was introduced in~\cite{CS:1991} (for arbitrary field coefficients). A problem with the ambiguity of mapping cones related with the definition given in~\cite{CS:1991} was improved in~\cite{You:1997} (see also \cite{BSY:2010}). In ~\cite{FPS:2023}, a more flexible definition of cobordism group has been introduced.  We refer to~\cite[Section 1.1]{FPS:2023} for the definition and properties of cobordism group (see also Conventions in Section \ref{intro}). 

 In this note, a {\em self-dual} $\R$-complex will be a pair $(\mathcal F,S)$, where $\mathcal F$ is a bounded complex with constructible cohomology sheaves on $X$, and $S:\mathcal F\otimes\mathcal F\to \mathbb D_X$ is a perfect pairing, where $\mathbb D_X$ is the dualizing complex on $X$. Perfect pairing means that the induced homomorphism $\alpha_S:\mathcal F\to \mathcal D_X(\mathcal F)$ is an isomorphism, where $\mathcal D_X(-)=Hom(-,\mathbb D_X)$ is Verdier's duality functor on $X$. Hence, this definition is equivalent to the classical  notion of self-duality where a self-dual $\R$-complex is a pair $(\mathcal F,\alpha)$, where $\alpha$ is an isomorphism from $\mathcal F$ to $\mathcal D_X(\mathcal F)$.  
  
Theorem \ref{conj:BSY} and Theorem \ref{main2} not only hold for the cobordism group introduced in \cite{FPS:2023}, but also for the cobordism group defined in~\cite{CS:1991},~\cite{You:1997}, \cite{BSY:2010}, as the properties we use in our proofs apply to both definitions of the cobordism group (see also Remark~\ref{re:all}).

\subsection{Polarizations at stalks of local systems}

Let $\varepsilon:Z\to X$ be a projective morphism of complex algebraic varieties, with $Z$ smooth of dimension $d$. Since $R\varepsilon_*$ commutes with Verdier duality for proper maps (see \cite[Chapter III]{KS:1994}), we have that $R\varepsilon_*\R_Z[d]$ inherits a self-duality
\begin{equation}
\label{eq:self1}
S:=\mathrm{tr}\circ R\varepsilon_*\sigma_Z:R\varepsilon_*\R_Z[d]\otimes R\varepsilon_*\R_Z[d]\to\mathbb D_X,
\end{equation}
where $\mathrm{tr}: R\varepsilon_*\mathbb D_Z=\varepsilon_{!}\varepsilon^!\mathbb D_X\to\mathbb D_X$ is the trace morphism defined by adjunction.

In this section the complexes $C$ appearing will be direct sums of intersection cohomology complexes associated with local systems. Given such a complex $C$ and a subvariety $Y_j\subset X$, we denote by $C_{Y_j}$ the direct sum of those direct summands of $C$ whose support is exactly $Y_j$.

Theorem \ref{teo:DecTheo} gives the direct sum decomposition 
\begin{equation}
\label{eq:decom0}
R\varepsilon_*\R_Z[d]\cong \bigoplus_{i=-M}^M {^p\hH}^i(R\varepsilon_*\R_{Z}[d])[-i].
\end{equation}
The self-duality~(\ref{eq:self1}) induces a self-duality
\begin{equation}
\label{eq:self2}
{^p\hH}^0(S):{^p\hH}^0(R\varepsilon_*\R_{Z}[d])\otimes  {^p\hH}^0(R\varepsilon_*\R_{Z}[d])    \to \mathbb D_X,
\end{equation}
and by \cite{You:1997} (compare with \cite[Lemma 3.3]{CS:1991}, with~\cite[Proposition 2.14]{SW:2019} for a similar result, and with~\cite{FPS:2023} for a proof for the new definition of cobordism group), the following equality of cobordism classes holds:
\begin{equation}
\label{eq:cobor1}
 [(R\varepsilon_*\R_Z[d],S)] =[({^p\hH}^0(R\varepsilon_*\R_{Z}[d]),{^p\hH}^0(S))]\in \Omega_\R(X).
\end{equation}

Let $\eta$ denote the first Chern class of a relative ample bundle for $\varepsilon$. By~\cite{BBD:1982} (see~\cite{Saito:1989a} and \cite{DeCM:2005} for characteristic $0$ proofs, the later only using classical Hodge theory) Hard-Lefschetz Theorem is satisfied for the decomposition above. That is, $\eta$ induces isomorphisms 
\begin{equation}
\label{eq:HLiso}
\eta^i:{^p\hH}^{-i}(R\varepsilon_*\R_{Z}[d])\to {^p\hH}^i(R\varepsilon_*\R_{Z}[d]),
\end{equation}
and we have the direct sum decomposition 
\begin{equation}
\label{eq:HLdecomp}
{^p\hH}^{-i}(R\varepsilon_*\R_{Z}[d])\cong \bigoplus_{l\geq 0} \eta^{l}\mathcal{P}^{-i-2l}(R\varepsilon_*\R_{Z}[d])
\end{equation}
for every non-negative $i$, where $\mathcal{P}^{-i-2l}(R\varepsilon_*\R_{Z}[d])$ denotes the primitive part of ${^p\hH}^{-i-2l}(R\varepsilon_*\R_{Z}[d])$. We remind the reader that $\mathcal{P}^{-i}(R\varepsilon_*\R_{Z}[d])$ is defined to be the kernel of $\eta^{i+1}:{^p\hH}^{-i}(R\varepsilon_*\R_{Z}[d])\to {^p\hH}^{i+2}(R\varepsilon_*\R_{Z}[d])$ in the abelian category of perverse sheaves. 
For $i=0$ this decomposition is orthogonal for the self-duality~(\ref{eq:self2}).

Theorem \ref{teo:DecTheo} also implies that each $\mathcal{P}^{-i}(R\varepsilon_*\R_{Z}[d])$ is a direct sum of simple intersection cohomology complexes, so, we have the decomposition 
\begin{equation}
\label{eq:decomsupp}
\mathcal{P}^{-i}(R\varepsilon_*\R_{Z}[d])\cong\bigoplus_{j\in J}\mathcal{P}^{-i}(R\varepsilon_*\R_{Z}[d])_{Y_j},
\end{equation}
where $\{Y_j\}_{j\in J}$ is the collection of possible supports. Since any morphism between two simple perverse sheaves with different strict support vanishes, this decomposition is also orthogonal for the self-duality~(\ref{eq:self2}). As a consequence we get the following equality in $\Omega_{\R}(X)$:
\begin{equation}
\label{eq:cobor2}
[({^p\hH}^0(R\varepsilon_*\R_{Z}[d]),{^p\hH}^0(S))]=\sum_{j\in J}\sum_{l\geq 0} [(\eta^l\mathcal{P}^{-2l}(R\varepsilon_*\R_{Z}[d])_{Y_j},S^l_{Y_j})],
\end{equation}
where $S^l_{Y_j}$ denotes the restriction of the perfect pairing ${^p\hH}^0(S)$ to the orthogonal direct summand $\eta^l\mathcal{P}^{-2l}(R\varepsilon_*\R_{Z}[d])_{Y_j}$.

The complex $\mathcal{P}^{-2l}(R\varepsilon_*\R_{Z}[d])_{Y_j}$ is isomorphic to $IC_{Y_j}(\mathcal{L})$, for a certain local system $\mathcal{L}$ in a Zariski open subset $U_j$ of $Y_j$. For each $y\in U_j$ we consider the inclusion $i_y:\{y\}\to Y_j$. We have the identification $i_{y}^!IC_{Y_j}(\mathcal{L})[\dim Y_j]\cong  \mathcal L_y$ and, applying the functor $i_y^!$, we obtain a perfect pairing 
$$Q_{\alpha_{Y_j},y}:=i_y^!S^l_{Y_j}:\mathcal L_y\otimes \mathcal L_y\to \R.$$ 


\begin{lem}
\label{lem:polarization}
Define $d(j,i):=\dim Z-\dim Y_j-i$. For any point $y\in U_j$ the stalk $\mathcal{L}_y$ has a natural structure of pure $\R$-Hodge structure of weight $d(j,2l)$ and $$(-1)^{(1/2)d(j,2l)(d(j,2l)+1)+(d-\dim Y_j)\dim Y_j}Q_{\alpha_{Y_j},y}$$
is a polarization.
\end{lem}

\begin{re}
 A proof of this lemma is also possible using M. Saito theory of Hodge modules as follows: By~\cite[Theorem 5.3.1]{Saito:1989a} $IC_{Y_j}(\mathcal{L})$ underlies a polarized pure Hodge module, whose polarization is the perfect pairing $S^l_{Y_i}$ up to a sign which is precisely determined. Such a pure Hodge module corresponds to a polarized variation of Hodge structures whose local system is $\mathcal{L}$. A dictionary comparing the signs of polarizations of pure Hodge modules and polarizations of their corresponding variation of pure Hodge structures is provided in~\cite[5.2.12]{Saito:1989a}. Here, we have to notice that since in our convention for polarization we insert Weil's operator on the left (see Conventions in Section \ref{intro}), and in Saito's convention it is inserted in the right, one need to multiply by the extra sign $(-1)^w$, where $w$ is the weight of the variation of pure Hodge structures. Then, the sign dictionary is the following: if a perfect pairing $S$ induces a polarization of a variation of pure Hodge structures of weight $w$ and support of pure dimension $d$, then $(-1)^{(1/2)d(d-1)+w}S$ induces a polarization of the corresponding pure Hodge module. 
 
However, for our proof, no Hodge modules or variations of Hodge structures are really needed. For us it is enough to understand a single stalk of $\mathcal L$. So, below we prove the lemma using elementary computations based on classical Hodge theory. The spirit of our proof reminds of the way that classical Hodge theory is used in~\cite{DeCM:2005}.
\end{re}

\proof[Proof of Lemma~\ref{lem:polarization}]

First we prove the lemma for the case of the structure morphism $\varepsilon:Z\to X=\{pt\}$.  Due to the shift, the intersection form 
$$Q:H^d(Z,\R)\otimes H^d(Z,\R)\to H^{2d}(Z,\R)\cong\R$$     
that this pairing induces equals $(-1)^d \langle -,-\rangle$, where   
$$\langle -,-\rangle:H^d(Z,\R)\otimes H^d(Z,\R)\to H^{2d}(Z,\R)\cong\R$$
is the usual intersection form, induced by the pairing $\R_Z\otimes\R_Z\to\R_Z$. Indeed, let $A_Z^\bullet$ be the de-Rham complex of $Z$. We have a chain of isomorphisms
$$A_Z^\bullet[d]\otimes A_Z^\bullet[d]\cong A_Z^\bullet\otimes\R[d]\otimes A_Z^\bullet\otimes\R[d]\cong A_Z^\bullet\otimes  A_Z^\bullet\otimes\R[d]\otimes\R[d]\cong A_Z^\bullet\otimes A_Z^\bullet\otimes\R[2d],$$
of which the second maps $\beta\otimes \lambda[d]\otimes\gamma\otimes\mu[d]\to (-1)^{dl} \beta\otimes\gamma\otimes\lambda[d]\otimes\mu[d]$ for $\beta\otimes \lambda[d]\otimes\gamma\otimes\mu[d]\in A_Z^k\otimes\R[d]\otimes A_Z^l\otimes\R[d]$. This induces a $(-1)^{d^2}=(-1)^d$ sign comparing the pairings $Q$ and $\langle -,-\rangle$.

In this case the stalk $\mathcal L_y$ is identified with $\eta^l\mathcal P^{d-2l}(Z)$, where $\mathcal P^{d-2l}(Z)$ denotes the $\eta$-primitive part of the cohomology $H^{d-2l}(Z,\R)$. By the classical Hodge-Riemann bilinear relations (see for instance~\cite[Theorem 1.33]{PeSte:2008}), we have that $(-1)^{(1/2)(d-2l)(d-2l-1)}\langle -,-\rangle$ is a polarization of $\mathcal P^{d-2l}(Z)$, with the pure Hodge structure inherited from $H^{d-2l}(Z,\R)$. The lemma holds in this case because $(-1)^d=(-1)^{d-2l}$ and $(1/2)(d-2l)(d-2l-1)+d-2l=(1/2)(d-2l)(d-2l+1)$.

Now we consider the general case. Let $H_j\subset X$ be the intersection of $\dim Y_j$ generic hyperplane sections in $X$ and denote by $\iota:H_{j}\hookrightarrow X$ the inclusion. By genericity we have that 
\begin{enumerate}
\item $Z_{H_j}:=\varepsilon^{-1}(H_j)$ is smooth and $\varepsilon|_{Z_{H_j}}:Z_{H_j}\to H_j$ is a resolution of singularities.
\item Set $c:=\dim Y_j $. The intersection $Y_{j'}\cap H_j$ is of dimension $\dim Y_{j'}-c$ and empty if $c>\dim Y_{j'} $. If the intersection is not empty then $U_{j'}\cap H_j$ is dense in $Y_{j'}\cap H_j$.
\item There is a tubular neighborhood $T(H_j)$ in $X$ and a continuous retraction map $\pi:T(H_j)\to H_j$ such that $\pi$ is topologically equivalent to a real vector bundle over $H_j$ of rank $2c$, and such that for every $j'\in J$ such that $\dim Y_{j'}\geq c$ we have 
$Y_{j'}\cap T(H_j)=\pi^{-1}(Y_{j'}\cap H_j)$ and $U_{j'}\cap T(H_j)=\pi^{-1}(U_{j'}\cap H_j)$.
\end{enumerate}
By  Property (3) above and~\cite[5.4.1, 5.4.3]{GM2:1983}, we obtain that for every $j'$ such that $\dim Y_{j'} \geq c$ the complex $\iota^! \eta^{l}\mathcal{P}^{-2l}(R\varepsilon_*\R_{Z}[d])_{Y_{j'}}[\dim Y_j]$ is the intersection cohomology complex associated with the restriction to $U_{j'}\cap\ H_j$ of the local system corresponding to $\eta^{l}\mathcal{P}^{-2l}(R\varepsilon_*\R_{Z}[d])_{Y_{j'}}[\dim Y_{j}]$.

By~\cite[5.4.1, 1, Formula (13) of 1.13]{GM2:1983}, we have that $\iota^!\R_Z[d][\dim Y_j]=\R_{Z_{H_j}}[d-\dim Y_j]$ and 
$R(\varepsilon|_{Z_{H_j}})_*\R_{Z_{H_j}}[d-\dim Y_j]=\iota^! R\varepsilon_*\R_{Z}[d][\dim Y_j]$. Then, applying ${^p\hH}^{0}(\iota^! -)$ to the decompositions~(\ref{eq:decom0}), (\ref{eq:HLdecomp}) and (\ref{eq:decomsupp}), and noticing that $\iota^!$ transforms the shifted perverse sheaves appearing in the decomposition into shifted perverse sheaves, we obtain the decomposition 
$$
{^p\hH}^{0}(R(\varepsilon|_{Z_{H_j}})_*\R_{Z_{H_j}}[d-\dim Y_j])\cong \bigoplus_{l\geq 0}\bigoplus_{j'\in J} \iota^!\eta^{l}\mathcal{P}^{-2l}(R\varepsilon_*\R_{Z}[d])_{Y_{j'}}[\dim Y_j],$$
which is the $\eta$-primitive decomposition of ${^p\hH}^{0}(R(\varepsilon|_{Z_{H_j}})_*\R_{Z_{H_j}}[d-\dim Y_j])$.

We denote by $S^l_{Y_{j'}\cap H_j}$ the restriction of the perfect pairing
$R(\varepsilon|_{Z_{H_j}})_*S_{Z_{H_j}}$ to the orthogonal summand  $\iota^!\eta^{l}\mathcal{P}^{-2l}(R\varepsilon_*\R_{Z}[d])_{Y_{j'}}[\dim Y_j]$. We have the equality of perfect pairings
\begin{equation}
\label{eq:selfindent}
S^l_{Y_{j'}\cap H_j}=(-1)^{(d-\dim Y_j)\dim Y_j}\iota^!S^l_{Y_{j'}}.
\end{equation}
Here the sign comes by a reason analogous to the sing comparing $Q$ and $\langle -,-\rangle$ above, since the dimension of $Z_{H_j}$ is $d-\dim Y_j$ and $\iota^!\R_Z[d][\dim Y_j]=\iota^*\R_Z[d][-\dim Y_j]$.
Hence, the bilinear form $Q_{\alpha_{Y_{j'}},y}$ coincides with the bilinear form $Q_{\alpha_{Y_{j'}\cap H_j},y}$ associated with $e_y^!S^l_{Y_{j'}\cap H_j}$ up to the sign $(-1)^{(d-\dim Y_j )\dim Y_j}$. This reduces the proof to the case in which $Y_j$ equals a point.

Assume that $Y_j=\{y\}$. The vector space $\eta^l\mathcal{P}^{-2l}(R\varepsilon_*\R_Z[d])_{Y_j}$ inherits a perfect pairing $S^l_{Y_j}$, which coincides with the bilinear form $Q_{\alpha_{Y_j},y}$ since $Y_j=y$. 
By \cite[Corollary 2.1.7, Theorem 2.1.8]{DeCM:2005}, $\mathcal{P}^{-2l}(R\varepsilon_*\R_Z[d])_{Y_j}$ is a $\R$-Hodge structure of weight $d-2l$ and that $Q_{\alpha_{Y_j},y}$ is a polarization up to a sign. Below we reduce a self contained argument proving this and determining the sign.

Consider the structure morphisms $f:Z\to \{pt\}$, $g:X\to \{pt\}$. Notice that $\mathcal{P}^{-2l}(R\varepsilon_*\R_{Z}[d])_y$ is a direct summand both $R\varepsilon_*\R_Z[d]$ supported at $\{y\}$ and of the primitive part $\mathcal P^{d-2l}(Z)$ of $H^{d-2l}(Z,\R)=R^{-2l}f_*\R_Z[d]$. Since we have proven the lemma for the case of the structure morphims, the perfect pairing $Rf_*\sigma_Z$ induces a bilinear form $Q'_y$ on $\eta^l\mathcal P^{d-2l}(Z)$ such that the form $(-1)^{(1/2)(d-2l)(d-2l+1)}Q'_y$ is a polarization. Then the proof of the lemma is finished because the restriction of $Q'_y$ to $\eta^l\mathcal{P}^{-2l}(R\varepsilon_*\R_{Z}[d])_y$ coincides with $Q_{\alpha_{Y_j},y}$. The last claim holds because we have the equality of perfect pairings
$$Rf_*\sigma_Z=Rg_*R\varepsilon_*\sigma_Z,$$ 
and since $\eta^l\mathcal{P}^{-2l}(R\varepsilon_*\R_{Z}[d])_y$ is supported at a point $y$ the functor $Rg_*$ restricted to it takes global sections and identifies the perfect pairing $Rf_*\sigma_Z$ with the the restriction of $R\varepsilon_*\sigma_Z$ to $\eta^l\mathcal{P}^{-2l}(R\varepsilon_*\R_{Z}[d])_y$.
\medskip

J. Sch\"urmann pointed out to us that an alternative proof should follow from \cite[Theorem 5.3.1,  Remark 5.3.12]{Saito:1989a}.
\endproof

\begin{re}
The use of hyperplane sections in the above lemma does not force projectivity assumptions. Indeed the statement is local in $X$, and $X$ can be covered by affine patches that can be completed to projective varieties for which the proof works (the completion is needed because the compactness of the resolution $Z$ is used in the last part of the proof). 
\end{re}

\subsection{Independence of the polarizing self-duality in the cobordism group}\label{sec:representation}
Although the local systems appearing in the previous section were all of them semi-simple, the results of this section hold without 
the semi-simplicity condition, and they (and their proof) may be useful in this generality elsewhere.

\begin{defi}
\label{def:polari}
Let $\mathcal V$ be a $\R$-perverse sheaf with strict support in an irreducible variety $X$. Let $U$ be a Zariski open subset such that $\mathcal V|_U=\mathcal L [\dim  X]$ where $\mathcal L$ is a $\R$-local system. Assume that for a given $x\in U$ the fibre $\mathcal L_x$ is endowed with a pure Hodge structure.  Let $\alpha:\mathcal{V}\otimes \mathcal{V}\to \mathbb D_{X}$ be a perfect paring. Then $\alpha$ is called a {\em polarizing self-duality for the Hodge structure at $x$} if $\alpha$ induces a polarization $Q_{\alpha,x}\colon\mathcal{L}_x\times \mathcal{L}_x\to\R$ of the Hodge structure. If the negative of the self-duality is polarizing, then we say that the self-duality is {\em $(-1)$-polarizing for the Hodge structure at $x$}.
\end{defi}

\begin{re}
\label{re:independence}
Set $d_j:=\dim Y_j$. Lemma~\ref{lem:polarization} states that the pairing $S^l_{Y_j}$ on $\mathcal{P}^{-2l}(R\varepsilon_*\R_{Z}[d])_{Y_j}$ induces a 
$(-1)^{(1/2)d(j,2l)(d(j,2l)+1)+(d-d_j)d_j}$-polarizing self-duality for the Hodge structure induced at $\mathcal L_y$ for any point $y\in U_i$. In other words, the fact that $S^l_{Y_j}$ is polarizing or $(-1)$-polarizing does not depend on the point $y\in U_j$. 
\end{re}

\begin{lem}
\label{lem:easydetermination}
Let $\alpha$ and $\alpha'$ be polarizing self-dualities of $\mathcal{V}$ for the same Hodge structure at $x$. Then $(\mathcal{V},\alpha)$ and $(\mathcal{V},\alpha')$ represent the same element in $\Omega_\R(X)$. 
\end{lem}
\proof
A different proof is provided in~\cite{FPS:2023}. The proof here is a based in a Representation Theory proposition that may be of independent interest. 

For any $s\in [0,1]$ the morphism $s\alpha+(1-s)\alpha'$ is a polarizing self-duality. Indeed, since both $\alpha$ and $\alpha'$ are polarizing there exits a common open subset $U$ of $X$ such that $Q_{\alpha,x}$ and $Q_{\alpha',x}$ are polarizations for any $x\in U$. Therefore an straightforward check of the conditions of~\cite[Definition 2.9]{PeSte:2008} imply that for any $x\in U$, the bilinear form $Q_{s\alpha+(1-s)\alpha',x}=sQ_{\alpha,x}+(1-s)Q_{\alpha',x}$ is a polarization of the Hodge structure $\mathcal{L}_x$. This implies that $Q_{s\alpha+(1-s)\alpha',x}$ is non-degenerate for any $s$. Then, in order to complete the proof it is enough to use Proposition~\ref{prop:rep} below.
\endproof

\begin{re}
\label{re:all}
Note that Lemma \ref{lem:easydetermination} holds for all the possible definitions of $\Omega_\R(X)$ since the cobordism relation is not used in the proof as was pointed out in the introduction. 
\end{re}


\begin{propo}
\label{prop:rep}
Let $\mathcal L$ be a $\R$-local system on $U$ and $Q_s$ a smooth $1$-parameter family of non-degenerate symmetric or anti-simmetric pairings of $\mathcal L$ parametrized by an interval $I$. Then $I$ admits an open cover $I=\cup_{j\in J} I_j$ such that for any $j$ there exists a smooth family of automorphisms $T(s)$ of $\mathcal L$ , $s\in I_j$, so that $Q_s(T(s)(-),T(s)(-))$ is independent of $s$. Consequently, for any $s,s'\in I$ there is an automorphism $T_{s,s'}$ of $\mathcal L$ such that $$Q_{s'}(T_{s,s'}(-),T_{s,s'}(-))=Q_s(-,-).$$
\end{propo}
\proof
The local system $\mathcal{L}$ is a representation $\rho:\pi_1(U,y)\to GL(\mathcal{L}_y)$ which is orthogonal for $Q_{s,y}$ for any $s$.

By Gram-Schmidt process in the symmetric case, and by the proof of uniqueness of non-degenerate anti-symmetric real bilinear forms, for any $s_0\in I$ there exists a neighborhood $I_{s_0}$ of $s_0$ in $I$ and a smooth family of automorphisms $N(s)$, $s\in I_{s_0}$ (not necessarily compatible with the monodromy) such that $N(s_0)=Id$ and $Q_{s,y}(N(s)(-),N(s)(-))=Q_{s_0,y}(-,-)$ for all $s$. Considering $M_{\gamma(s)}:=N(s)^{-1}\rho(\gamma)N(s)$ we obtain a smooth family of real orthogonal representations for $Q_{s_0,y}$. 

If $Q_{s_0,y}$ is symmetric and $(n,m)$ is its signature of $Q_{s_0,y}$ we define $W$ to be the diagonal matrix $I_{n,m}$ be the diagonal matrix of size $n+m$ such that the $(i,i)$ component equals $1$ if $i\leq n$ and $-1$ if $i>n$. If $Q_{s_0,y}$ is anti-symmetric the rank of the local system is even and we define $W$ to be the matrix 
\begin{equation}
\begin{pmatrix}
 0 & I\\
 -I & 0
\end{pmatrix}
\end{equation}
where $I$ denotes the identity matrix. Denote $O(W)$ the orthogonal group for the bilinear form $W$. Then, the proposition is reduced to the following claim:

\textsc{Claim}: Let $M_{\gamma}(s):\pi_1(U_j,y)\to O(W)$ be a smooth family of orthogonal representations for the quadratic form $W$. If there exists a family $N(s)$ of invertible matrices such that $N(0)=Id$ and we have the conjugation $M_{\gamma}(s)=(N(s))^{-1}M_\gamma(0) N(s)$ for any $\gamma\in \pi_1(U_j,y)$, then there exists a family $P(s)$ of orthogonal matrices such that $M_{\gamma}(s)=(P(s))^{-1}M_\gamma(0) P(s)$ for any $\gamma\in \pi_1(U_j,y)$.

In order to make the proof of the claim more self-contained (since it uses different ideas than the rest of the paper), we recall an elementary fact about Lie groups (if the reader wishes to consider only the case of matrix Lie groups, this is the case that will be used, but the discussion holds in general). The reader may refer to~\cite{Lee:2002} for generalities on smooth manifolds, and Lie groups; flows associated to time dependent vector fields are discussed in Exercise 12-7 of loc. cit. 

Let $G$ be a matrix Lie group. Its Lie algebra $\mathfrak{g}$ is identified with the spaces of left-invariant vector fields. Denote by $L_g:G\to G$ the left multiplication by $g$. We denote the set of smooth paths $g:[0,1]\to G$ by $C^\infty([0,1],G)$. A path $v(t)\in C^\infty([0,1],\mathfrak{g})$ is viewed as a left-invariant time dependent vector field, that is, a vector field in $G\times [0,1]$ such that its $[0,1]$-component is the unit vector field $\frac{\partial }{\partial t}$ in positive direction, and such that its $G$-component is invariant by the action of $G$ by left multiplication. 
We say that $v(t)$ is {\em integrable} if the integral flow associated with it is defined in the domain $G\times [0,1]$. In order to check integrability, by left invariance, it is enough to check the existence of an integral curve whose domain is $[0,1]$.  
We denote by $C^\infty([0,1],\mathfrak{g})^{int}$ the set of maps giving rise to integrable left invariant time dependent vector fields. 

We define the {\em left bijection}
$$\mathcal{L}:C^\infty([0,1],G)\to G\times C^\infty([0,1],\mathfrak{g})^{int}$$
as follows. Given a smooth path $g:[0,1]\to G$ we define $g':[0,1]\to\mathfrak{g}$ by the formula 
$$g'(s):=DL_{g(s)^{-1}}(g(s))(\frac{dg(u)}{du}|_{u=s}).$$ Given $g(s)\in C^\infty([0,1],G)$ we define $\mathcal{L}(g(s))$ to be the pair $(g(0),g'(s))$. Conversely, given a pair $(g_0,v(t))\in G\times C^\infty([0,1],\mathfrak{g})^{int}$, we view $v(t)$ as a left-invariant time dependent vector field and define $\mathcal{L}^{-1}(g_0,v(t))$ to be the unique integral curve of the time dependent vector field $v(t)$ with initial point $g(0)=g_0$. 

The set $C^\infty([0,1],G)$ has a group structure. Given $h\in G$ and $g(s)\in C^\infty([0,1],G)$, define $h(s):=g(s)^{-1}h g(s)$. A straightforward Lie group computation shows the formula
\begin{equation}
\label{eq:formulaLie}
h'(s)=g'(s)-h(s)^{-1}g'(s)h(s).
\end{equation}
Indeed, since left multiplcation by a matrix is a linear transformation at the space of matrices we may write $$DL_{g(s)^{-1}}(g(s))(\frac{dg(u)}{du}|_{u=s})=g(s)^{-1}(\frac{dg(u)}{du}|_{u=s}).$$ Using this, Leibnitz rule for derivation and the formula for the derivation of the inverse we obtain 
\begin{align*}
&h'(s)=h(s)^{-1}(\frac{dh(u)}{du}|_{u=s})=\\
&=h(s)^{-1}(-g(s)^{-1}(\frac{dg(u)}{du}|_{u=s})g(s)^{-1}hg(s)+g(s)^{-1}h(\frac{dg(u)}{du}|_{u=s}))=\\
&=-h(s)^{-1}g(s)^{-1}(\frac{dg(u)}{du}|_{u=s})g(s)^{-1}hg(s)+g(s)^{-1}h^{-1}g(s)g(s)^{-1}h(\frac{dg(u)}{du}|_{u=s})=\\
&=-h(s)^{-1}g'(s)h(s)+g'(s).
\end{align*}

Let $W$ be $(-1)^\beta$-symmetric. Then we have $W^2=(-1)^\beta Id$ and the the Lie algebra $\mathfrak{o}(W)$ of the Lie group is $O(W)$ the subspace of matrices $N$ satisfying $(-1)^{\beta+1}WN^tW=N$. So $\mathfrak{o}(W)$ is the eigenspace for eigenvalue $1$ of the involution $N\mapsto (-1)^{\beta+1}WN^tW$ in the real vector space of square matrices. The vector space of square matrices splits as the direct sum of the eigenspaces with eigenvalues $+1$ and $-1$ respectively for the involution. Given any square matrix $N$ we decompose it accordingly as 
$N=N_++N_-$.

Formula~(\ref{eq:formulaLie}) applied to $M_{\gamma}(s)=N(s)^{-1}M_\gamma(0) N(s)$ yields
$$
M_\gamma'(s)=N'(s)-M_\gamma(s)^{-1}N'(s)M_\gamma(s).
$$
We have 
\begin{equation}
\label{eq:matrixeq}
M_\gamma'(s)=M_\gamma'(s)_+=N'(s)_+ -M_\gamma(s)^{-1}N'(s)_+M_\gamma(s).
\end{equation}
The first equality is because $M_\gamma(s)$ belongs to $O(W)$ for all $s$. For the second equality write $N'(s)=N'(s)_++N'(s)_-$. It is enough to show the equalities 
$$(M_\gamma(s)^{-1}N'(s)M_\gamma(s))_+=M_\gamma(s)^{-1}N'(s)_+M_\gamma(s),$$
$$(M_\gamma(s)^{-1}N'(s)M_\gamma(s))_-=M_\gamma(s)^{-1}N'(s)_-M_\gamma(s),$$
but they follow from an elementary matrix computation using that since $M_\gamma(s)\in O(W)$ we have $M_\gamma(s)^{-1}=(-1)^\beta WM_\gamma(s)^tW$.

Define $P(s):=\mathcal{L}^{-1}(Id,N'(s)_+)$. Since $N'(s)_+\in\mathfrak{o}(W)$ we have $P(s)\in O(W)$. The equality $M_{\gamma}(s)=(P(s))^{-1}M_\gamma(0) P(s)$ is obtained applying $\mathcal{L}^{-1}(Id,\bullet)$ to Equation~(\ref{eq:matrixeq}).
\endproof

\subsection{An identity in the cobordism group}

In Lemma \ref{lem:rewriting} we use the following notation: Certain semi-simple perverse sheaves $\mathcal V$ with strict support  are endowed with polarizable Hodge structures at their stalks at generic points by Lemma~\ref{lem:polarization}. We denote by $[\mathcal V,\oplus]$ the class in $\Omega_\R(X)$ represented by $\mathcal V$ together with a polarizing self-duality (this definition makes sense by Remark~\ref{re:independence} and Lemma~\ref{lem:easydetermination}). We denote by $[\mathcal V,\ominus]$ the class with a $(-1)$-polarizing self-duality. If a semi-simple perverse sheaf $\mathcal W$ together with a self-duality is a direct sum of polarizing semi-simple perverse sheaves with strict support, then we denote by $[\mathcal W,\oplus]$ its class in $\Omega_\R(X)$; we denote by $[\mathcal W,\ominus]$ the class of the opposite self-duality.

\begin{lem}
\label{lem:rewriting}
Define $\beta_{d,j}:=(d-\dim Y_j)\dim Y_j$. We have the following equality in $\Omega_\R(X)$:
\begin{equation}
\label{eq:cobor3}
\begin{split}
 [({^p\hH}^0(R\varepsilon_*\R_{Z}[d]),{^p\hH}^0(S))]= \\
 =\sum_{j\in J}\sum_{i=-M}^M (-1)^{(1/2)d(j,-i)(d(j,-i)+1)+\beta_{d,j}} [{^p\hH}^i(R\varepsilon_*\R_{Z}[d])_{Y_j},\oplus]=\\
 =\sum_{j\in J} (\sum_{i=\text{even}}(-1)^{(1/2)d(j,-i)(d(j,-i)+1)+\beta_{d,j}} [{^p\hH}^i(R\varepsilon_*\R_{Z}[d])_{Y_j},\oplus]+\\
 +\sum_{i=\text{odd}}(-1)^{(1/2)d(j,-i)(d(j,-i)+1)+\beta_{d,j}} [{^p\hH}^i(R\varepsilon_*\R_{Z}[d])_{Y_j},\ominus]).
\end{split}
\end{equation}
\end{lem}
\proof

For each $j\in J$, Equation (\ref{eq:HLdecomp}) gives the decomposition \[{^p\hH}^0(R\varepsilon_*\R_{Z}[d])_{Y_j}=\bigoplus_{l\geq 0} \eta^l\mathcal P^{-2l}(R\varepsilon_*\R_{Z}[d])_{Y_j}\] where the decomposition is orthogonal for the self-duality of  ${^p\hH}^0(R\varepsilon_*\R_{Z}[d])_{Y_j}$ and the     self-duality of $\eta^l\mathcal P^{-2l}(R\varepsilon_*\R_{Z}[d])_{Y_j}$, induced from the self-duality of  ${^p\hH}^0(R\varepsilon_*\R_{Z}[d])_{Y_j}$ is $(-1)^{(1/2)d(j,2l)(d(j,2l)+1)+\beta_{d,j}}$-polarizing by Lemma~\ref{lem:polarization}. Then we have the equality
\begin{equation}
\label{eq:qwerty}
\begin{split}
[{^p\hH}^0(R\varepsilon_*\R_{Z}[d])_{Y_j}]=
 \sum_{l\geq 0}(-1)^{(1/2)d(j,2l)(d(j,2l)+1)+\beta_{d,j}}[\mathcal P^{-2l}(R\varepsilon_*\R_{Z}[d])_{Y_j},\oplus]\in \Omega_\R(X).
 \end{split}
\end{equation} 

By Equations~(\ref{eq:HLiso}) and~(\ref{eq:HLdecomp}), 
$$[{^p\hH}^{-i}(R\varepsilon_*\R_{Z}[d])_{Y_j},\oplus]=[{^p\hH}^{i}(R\varepsilon_*\R_{Z}[d])_{Y_j},\oplus]=\sum_{l\geq 0} [\mathcal{P}^{-i-2l}(R\varepsilon_*\R_{Z}[d])_{Y_j},\oplus]$$
in $\Omega_\R(X)$ for any $i\geq 0$. 

Plugging this equalities into the middle term of~(\ref{eq:cobor3}), making the needed {\em cancellations} and comparing with the left hand side of~(\ref{eq:cobor3}) expressed as in~(\ref{eq:qwerty}), the first equality of~(\ref{eq:cobor3}) follows. 

We show here part of the cancellation process for the reader convenience. We write ${^p\hH}^{i}$ for $[{^p\hH}^{i}(R\varepsilon_*\R_{Z}[d])_{Y_j},\oplus]$ and $\mathcal P^{-2l}$ for $[\mathcal P^{-2l}(R\varepsilon_*\R_{Z}[d])_{Y_j},\oplus]$. Then, using Hard-Lefschetz theorem, the following equalities hold
\begin{align*}
&(-1)^{(1/2)d(j,0)(d(j,0)+1)+\beta_{d,j}}\sum_{i=-M}^M (-1)^{(1/2)d(j,-2i)(d(j,-2i)+1)+\beta_{d,j}} {^p\hH}^{2i}=\\
&={^p\hH}^0-2{^p\hH}^{-2}+2{^p\hH}^{-4}-2{^p\hH}^{-6}+...=\\
&=(\mathcal P^0+\mathcal P^{-2}+\mathcal P^{-4}+\mathcal P^{-6}+...)-2(\mathcal P^{-2}+\mathcal P^{-4}+\mathcal P^{-6}+...)+\\
&+2(\mathcal P^{-4}+\mathcal P^{-6}+...)-2
(\mathcal P^{-6}+...)+...=\\
&=\mathcal P^0-\big{(}(\mathcal P^{-2}+\mathcal P^{-4}+\mathcal P^{-6}+...)-2(\mathcal P^{-4}+\mathcal P^{-6}+...)+\\
&+2(\mathcal P^{-6}+...)-...\big{)}=\mathcal P^0-\mathcal P^{-2}+\mathcal P^{-4}-\mathcal P^{-6}+...=\\
&=(-1)^{(1/2)d(j,0)(d(j,0)+1)+\beta_{d,j}}\sum_{l\geq 0}(-1)^{(1/2)d(j,2l)(d(j,2l)+1)+\beta_{d,j}}\mathcal P^{-2l}.
\end{align*}
This shows the cancellations in the middle term of formula~(\ref{eq:cobor3}) when $i$ is even. The $i$ odd part of the middle term of formula~(\ref{eq:cobor3}) cancels completely by a similar process, and by this complete cancellation the equality with the right hand side follows.
\endproof

\subsection{Proof of Theorem~\ref{main2}}\label{sec:proof}
For any $k$ the variety $X_k$ is a disjoint union of smooth varieties of different dimensions. By $d_k$ we denote the function that assigns to each connected component of $X_k$ its dimension, and given a complex of sheaves $C$ on $X_k$ we denote by $C[d_k]$ the same complex, shifted at the dimension in each connected component.

We have to prove that $sd_\R([X])-[IC_X]=0$ in $\Omega_\R(X)$. Equation~(\ref{Eq1}) implies
\begin{equation}
\label{eq:1111}
sd_\R([X])=sd_\R([\tilde X])+\sum_{i=1}^n sd_\R([X_{0,i}])+\sum_{k=1}^n (-1)^k sd_\R([X_k]).
\end{equation}
Since at the cobordism group $\Omega_\R(X)$ only the $0$-th perverse cohomology matters, by \cite{You:1997} (see also \cite[Lemma 3.3]{CS:1991}, \cite{FPS:2023}), then
$$sd_\R([\tilde X])+\sum_{i=1}^n sd_\R([X_{0,i}])=[IC_X]+\sum_{j=1}^N[R\varepsilon_*\R_{X_0}[d_0]_{\Sigma_j}]=$$
$$=[IC_X]+\sum_{j=1}^N[{^p\hH}^0(R\varepsilon_*\R_{X_0}[d_0])_{\Sigma_j}]$$
and
$$sd_\R([X_k])=\sum_{j=1}^N[R\varepsilon_*\R_{X_k}[d_k]_{\Sigma_j}]=\sum_{j=1}^N[{^p\hH}^0(R\varepsilon_*\R_{X_k}[d_k])_{\Sigma_j}]$$
for every $k>0$. 

Substituting the above expressions and using Lemma~\ref{lem:rewriting}, by setting $d_j:=\dim \Sigma_j$, we obtain 
$$sd_\R([X])-[IC_X]=$$

\begin{align*}
&=\sum_{j=1}^N \sum_{k=0}^n \sum_{i=-M_k}^{M_k}(-1)^{k+(1/2)(d_k-d_j+i)(d_k-d_j+i+1)+\beta_{d_k,j}} [{^p\hH}^i(R\varepsilon_*\R_{X_k}[d_k])_{\Sigma_j},\oplus]=\\
&=\sum_{j=1}^N \sum_{k=0}^n \sum_{q=0}^{2d_k}(-1)^{k+(1/2)(q-d_j)(q-d_j+1)+\beta_{d_k,j}} [{^p\hH}^q(R\varepsilon_*\R_{X_k})_{\Sigma_j},\oplus]=\\
&=\sum_{j=1}^N \sum_{k=0}^n \sum_{q-d_k=\text{even}}(-1)^{k+(1/2)(q-d_j)(q-d_j+1)+\beta_{d_k,j}} [{^p\hH}^q(R\varepsilon_*\R_{X_k})_{\Sigma_j},\oplus]+\\
&+\sum_{j=1}^N \sum_{k=0}^n \sum_{q-d_k=\text{odd}}(-1)^{k+(1/2)(q-d_j)(q-d_j+1)+(d_k-d_j+1)d_j} [{^p\hH}^q(R\varepsilon_*\R_{X_k})_{\Sigma_j},\oplus].
\end{align*}
The last equality uses the last equality of Equation~(\ref{eq:cobor3}) when $d_j$ is odd.

The proof concludes noticing that for any $q$, $k,k'$ such that $d_k-q$ is even and $d_{k'}-q$ is odd the signs $(-1)^{\beta_{d_k,j}}$ and $(-1)^{(d_k-d_j+1)d_j}$ coincide. Then, since the sign $(-1)^{(1/2)(q-d_j)(q-d_j+1)}$ is constant for $q$ and $\Sigma_j$ fixed, we apply Theorem~\ref{teo:equivalentRHM} for each support $\Sigma_j$ and any $q$ and conclude using Lemma~\ref{lem:easydetermination}.

\end{document}